\documentclass[12pt]{amsart}
\usepackage{fullpage,url}
\usepackage[dvips]{hyperref}
 \usepackage{amscd}   
\usepackage{amsmath}

\DeclareFontEncoding{OT2}{}{} 

\newcommand{\pp}{{\mathfrak p}}
\newcommand{\bet}{{\mathfrak P}}

\DeclareMathOperator{\qh}{QsM_{K^{\sep}}}

\DeclareMathOperator{\lang}{\text{SCF}_{p,\nu}}
\DeclareMathOperator{\Coker}{Coker}
\DeclareMathOperator{\qua}{QsE_{K^{\sep}}}

\DeclareMathOperator{\Ker}{Ker}

\DeclareMathOperator{\sep}{sep}

\DeclareMathOperator{\tor}{tor}
\DeclareMathOperator{\alg}{alg}

\DeclareMathOperator{\trdeg}{tr deg}

\DeclareMathOperator{\End}{End_{K^{\sep}}}

\DeclareMathOperator{\Aut}{Aut_{K^{\sep}}}
\DeclareMathOperator{\Gal}{Gal}

\DeclareMathOperator{\ord}{ord}

\DeclareMathOperator{\Frac}{Frac}

\newcommand{\isomto}{\overset{\sim}{\rightarrow}}

\newtheorem{theorem}{Theorem}[section]
\newtheorem{lemma}[theorem]{Lemma}
\newtheorem{corollary}[theorem]{Corollary}
\newtheorem{proposition}[theorem]{Proposition}

\theoremstyle{definition}
\newtheorem{definition}[theorem]{Definition}
\newtheorem{notation}[theorem]{Notation}
\newtheorem{statement}[theorem]{Statement}

\newtheorem{example}[theorem]{Example}

\theoremstyle{remark}
\newtheorem{remark}[theorem]{Remark}

\newtheorem{sublemma}[theorem]{Sublemma}

\title{The Mordell-Lang Theorem for Drinfeld modules}
\author{Dragos Ghioca}
\begin{document}
\begin{abstract}
We study the quasi-endomorphism ring of certain infinitely definable subgroups 
in separably closed fields. Based on the results we obtain, we are able to prove 
a Mordell-Lang theorem for Drinfeld modules of finite characteristic. Using 
specialization arguments we prove also a Mordell-Lang theorem for Drinfeld 
modules of generic characteristic. 
\end{abstract}
\maketitle
\section{Introduction}
Faltings proved the Mordell-Lang Conjecture in the following form (see 
\cite{Fal}).
\begin{theorem}[Faltings]
\label{T:F}
Let $G$ be a semiabelian variety defined over the field of complex numbers 
$\mathbb{C}$. Let $X\subset G$ be a closed subvariety and $\Gamma\subset 
G(\mathbb{C})$ a finitely generated subgroup of the group of $\mathbb{C}$-points 
on $G$. Then $X(\mathbb{C})\cap\Gamma$ is a finite union of cosets of subgroups 
of $\Gamma$.
\end{theorem}

If we try to formulate the Mordell-Lang Conjecture in the context of algebraic 
subvarieties contained in a power of the additive group scheme $\mathbb{G}_a$, 
the conclusion is either false (in the characteristic $0$ case, as shown by the 
curve $y=x^2$ which has an infinite intersection with the finitely generated 
subgroup $\mathbb{Z}\times\mathbb{Z}$, without being itself an additive 
algebraic group) or it is trivially true (in the characteristic $p>0$ case, 
because every finitely generated subgroup of a power of $\mathbb{G}_a$ is 
finite). In the third section we will present a nontrivial formulation of the 
Mordell-Lang conjecture for a power of the additive group in characteristic $p$ 
in the context of Drinfeld modules. We will replace the \emph{finitely generated 
subgroup} from the usual Mordell-Lang statement with a \emph{finitely generated 
$\phi$-submodule}, where $\phi$ is a Drinfeld module. We also strengthen the 
conclusion of the Mordell-Lang statement in our setting by asking that the 
\emph{subgroups} whose cosets are contained in the intersection of the algebraic 
variety with the finitely generated $\phi$-submodule are actually 
\emph{$\phi$-submodules}. 

In order to obtain the results of the present paper we need first to analyze the 
quasi-endomorphisms ring of a certain infinitely definable subgroup in the 
theory of separably closed fields. In this section we introduce the basic 
notation and results, while in the second section we prove the main result 
(Theorem ~\ref{T:t}) needed for the proof of Theorem ~\ref{T:DG} (the 
Mordell-Lang Theorem for Drinfeld modules of finite characteristic). Using 
specialization arguments we also prove a Mordell-Lang statement for Drinfeld 
modules of generic characteristic (Theorem ~\ref{T:T3}).

Everywhere in this paper, for two sets $A$ and $B$, $A\subset B$ means that $A$ 
is a subset, not neccessarily proper, of $B$.

Let $K$ be a finitely generated field of characteristic $p>0$. Let $\tau_0$ be 
the usual Frobenius, i.e. $\tau_0(x)=x^p$, for every $x$. We let $K\{\tau_0\}$ 
be the ring of all polynomials in $\tau_0$ with coefficients from $K$. If 
$f,g\in K\{\tau_0\}$, $fg$ will represent the composition of $f$ and $g$.

Fix an algebraic closure $K^{\alg}$ of $K$. Let $K^{\sep}$ be the separable 
closure of $K$ inside $K^{\alg}$. Let $\mathbb{F}_p^{\alg}$ be the algebraic 
closure of $\mathbb{F}_p$ inside $K^{\sep}$. 

Assume $\trdeg_{\mathbb{F}_p}K=\nu\ge 1$. Then $[K:K^p]=[K^{\sep}:K^{\sep 
^p}]=p^{\nu}>1$. The number $\nu$ is called the Ersov invariant of $K$. 

\begin{notation}
\label{N:tree}
Let $k$ be a positive integer. We denote by $p^{(k)}$ the set of functions 
$$f:\{1,\dots,k\}\rightarrow\{0,\dots,p-1\}.$$
\end{notation}

\begin{definition}
\label{D:p-basis}
A subset $B=\{b_1,\dots,b_{\nu}\}\subset K$ is called a $p$-basis of $K$, or 
equivalently, of $K^{\sep}$, if the following set of monomials,
$$\left\{ m_i=\prod_{j=1}^{\nu}b_j^{i(j)}\mid i\in p^{(\nu)}\right\}$$
forms a basis for $K/K^p$, or equivalently for $K^{\sep}/K^{\sep ^p}$.
\end{definition}
For the rest of this paper we fix a $p$-basis $B$ for $K$. For each $i\in 
p^{(\nu)}$, we denote by $\lambda_i:K^{\sep}\rightarrow K^{\sep}$ the unique 
functions with the property that for every $x\in K^{\sep}$,
$$x=\sum_{i\in p^{(\nu)}}\lambda_i(x)^p m_i.$$
We call these functions $\lambda_i$ the $\lambda$-functions of level $1$. For 
every $k\ge 2$ and for every choice of $i_1,\dots,i_k\in p^{(\nu)}$, 
$$\lambda_{i_1,i_2,\dots,i_k}=\lambda_{i_1}\circ\lambda_{i_2}\circ\dots\circ 
\lambda_{i_k}$$
is called a $\lambda$-function of level $k$.

\begin{definition}
\label{D:language}
We let $\lang$ be the theory of separably closed fields of characteristic $p$ 
and Ersov invariant $\nu$ in the language $$\mathfrak L _{p,\nu} 
=\{0,1,+,-,\cdot\}\cup\{b_1,\dots,b_{\nu}\}\cup \{\lambda_i\mid i\in 
p^{(\nu)}\}.$$ 
\end{definition}

From now on we consider a finitely generated field $K$ such that $K^{\sep}$ is a 
model of $\lang$. We let $L$ be an 
$\aleph_1$-saturated elementary extension of $K^{\sep}$. We are interested in 
studying infinitely definable subgroups $G$ of $(L,+)$, i.e. $G$ is possibly an infinite intersection of definable 
subgroups of $(L,+)$. If $k\ge 1$ and $G$ is an infinitely definable subgroup of 
$(L,+)$, then the \emph{relatively} definable subsets of $G^k$ (the cartesian 
product of $G$ with itself $k$ times) are the intersection of $G^k$ with 
definable subsets of $(L,+)^k$. If there is no risk of ambiguity, we will say a 
definable subset of $G^k$, instead of a relatively definable subset of $G^k$. 
The structure \emph{induced} by $L$ on $G$ over a set $S$ of parameters, is the 
set $G$ together with all the relatively $S$-definable subsets of the cartesian 
powers of $G$. We will consider only the case when the set $S$ of parameters 
equals $K^{\sep}$. Thus, when we will say a definable subset, we will mean a 
$K^{\sep}$-definable subset. Also, we will call \emph{additive} the subgroups of 
$(L,+)$.

\begin{definition}
\label{D:connected component}
For every infinitely definable subgroup $G$, the connected component of $G$, 
denoted $G^0$, is the intersection of all definable subgroups of finite index in 
$G$.
\end{definition}

\begin{definition}
\label{D:connected groups}
The group $G$ is connected if $G=G^0$.
\end{definition}

\begin{definition}
Let $G$ be an infinitely definable additive subgroup of $L$. We denote by 
$\End(G)$ the set of $K^{\sep}$-definable endomorphisms $f$ of $G$. If $G$ is a 
connected group, then the graph of $f$ is a connected subgroup of $G\times G$. 
The endomorphisms $f\in \End(G)$ that are both injective and surjective, 
form the group of $K^{\sep}$-automorphisms of $G$, denoted $\Aut(G)$.
\end{definition}
From now on, when we will say endomorphisms of $G$, we will refer to the 
elements of $\End(G)$ and when we will say automorphisms of $G$, we will refer 
to the elements of $\Aut(G)$.
\begin{definition}
\label{D:quasi-endomorphisms}
Let $G$ and $H$ be infinitely definable connected groups. We call the subgroup 
$\psi\subset G\times H$ a 
$K^{\sep}$-quasi-morphism from $G$ to $H$ if the following three properties are 
satisfied

1) $\psi$ is a connected, $K^{\sep}$-definable subgroup of $G\times H$.

2) the first projection $\pi_1(\psi)$ equals $G$.

3) the \emph{cokernel} $\Coker(\psi)=\{x\in H\mid (0,x)\in\psi\}$ is finite.

The set of all $K^{\sep}$-quasi-morphisms from $G$ to $H$ is denoted by 
$\qh(G,H)$. 

When $G=H$, we call $\psi$ a $K^{\sep}$-quasi-endomorphism of $G$. The set of 
all $K^{\sep}$-quasi-endomorphisms of $G$ is denoted by $\qua(G)$.
\end{definition}
For every infinitely definable additive subgroup $G$, when we will say 
quasi-endomorphisms of $G$, we will refer to the elements of $\qua(G)$. 

Let $f$ be an endomorphism of the connected group $G$. We interpret $f$ as a 
quasi-endomorphism of 
$G$ by $$f=\{(x,f(x))\mid x\in G\}\in\qua(G).$$

\begin{definition}
\label{D:qua is a ring}
Let $G$ be an infinitely definable connected group. We define the following two 
operations that will induce a ring structure on $\qua (G)$.

1) \emph{Addition.} For every $\psi_1,\psi_2\in\qua(G)$, we let $\psi_1+\psi_2$ 
be the connected component of the group 
$$\{(x,y)\in G\times G\mid\exists y_1,y_2\in G\text{ such that 
}(x,y_1)\in\psi_1\text{ and }(x,y_2)\in\psi_2\text{ and }y_1+y_2=y\}.$$

2) \emph{Composition.} For every $\psi_1,\psi_2\in\qua(G)$, we let 
$\psi_1\psi_2$ be the connected component of the group
$$\{(x,y)\in G\times G\mid\text{ there exists }z\in G\text{ such that }(x,z)\in 
\psi_2\text{ and }(z,y)\in\psi_1\}.$$
\end{definition}

\begin{definition}
\label{D:minimal group}
Let $G$ be an infinitely definable additive subgroup. Then $G$ is $c$-minimal if 
it is infinite and every definable subgroup of $G$ is either finite or has 
finite index.
\end{definition}

\begin{lemma}
\label{L:minimal}
If $G$ is a $c$-minimal connected group, then for all $f\in \End(G)\setminus\{0\}$, $f(G)=G$.
\end{lemma}

\begin{proof}
Because $f\in\End(G)$ and $G$ is connected, $f(G)$ is a definable, connected subgroup of $G$. thus either $f(G)=\{0\}$ or $f(G)$ is infinite. So, 
because $f\ne 0$, $f(G)$ is not finite. Then, because $G$ is $c$-minimal, $f(G)$ 
has finite index in $G$. Because $G$ is connected, we conclude that $f$ is surjective.
\end{proof}

The next result is proved in a larger generality in Section $4.4$ of 
\cite{Wagner}. Because for the case we are interested in we can give a simpler 
proof, we present our argument below.
\begin{proposition}
\label{P:quadivisionring}
If $G$ is a $c$-minimal, connected additive group, then $\qua(G)$ is a division 
ring.
\end{proposition}

\begin{proof}
Let $\psi\in\qua(G)\setminus\{0\}$. Let $\pi_2(\psi)$ be the projection of 
$\psi\subset G\times G$ on the second component. Then $\pi_2(\psi)$ is a 
definable subgroup of $G$. Because $\psi$ is connected and $\psi\ne 0$, 
$\pi_2(\psi)$ is not finite. Then, because $G$ is a $c$-minimal, connected 
group, $\pi_2(\psi)=G$.

Because $\pi_2(\psi)=G$ and $G$ is $c$-minimal, for each $y\in G$, the set
\begin{equation}
\label{E:preimage}
\{x\in G\mid (x,y)\in\psi\}
\end{equation}
is finite. We define $\phi=\{(y,x)\in G\times G\mid (x,y)\in\psi\}$. Because 
$\psi$ is a connected, $K^{\sep}$-definable subgroup of $G\times G$, then also 
$\phi$ is a connected, $K^{\sep}$-definable subgroup of $G\times G$. By 
construction, $\pi_1(\phi)=\pi_2(\psi)=G$. Let
$$\Coker(\phi)=\{x\in G\mid (0,x)\in \phi\}.$$
By construction of $\phi$, $\Coker(\phi)=\{x\in G\mid (x,0)\in\psi\}$. Using 
\eqref{E:preimage} for $y=0$, we conclude that $\Coker(\phi)$ is finite. Thus 
$\phi\in\qua(G)$. By definition of $\phi$, $\psi\phi$ (as defined in Definition 
~\ref{D:qua is a ring}) is the identity function on $G$. Thus $\qua(G)$ is a 
division ring.
\end{proof}

\begin{definition}
\label{D:fsharp}
Let $f\in K\{\tau_0\}\tau_0\setminus\{0\}$. We define 
$f^{\sharp}=f^{\sharp}(L)=\cap_{n\ge 
1}f^n(L)$.
\end{definition}

In \cite{Sca} and \cite{Blo} (Lemma $4.23$) it is proved the following result.
\begin{theorem}
\label{T:fsharpminimal}
If $f\in K\{\tau_0\}\tau_0\setminus\{0\}$, then $f^{\sharp}$ is $c$-minimal. In 
particular, 
$f^{\sharp}$ is infinite.
\end{theorem}

In \cite{Poi} it is proved the following result.
\begin{theorem}
\label{T:Ksepconnected}
The group $K^{\sep}$ is connected.
\end{theorem}

Because the image of a connected group through a definable map is also 
connected, we get the following result.
\begin{corollary}
\label{C:fksepconnected}
For every $f\in K\{\tau_0\}$, $f(K^{\sep})$ is connected.
\end{corollary}

\begin{lemma}
\label{L:fsharp connected}
If $f\in K\{\tau_0\}\tau_0\setminus\{0\}$, then $f^{\sharp}$ is connected.
\end{lemma}

\begin{proof}
It suffices to show that for every definable additive subgroup $G$ of $L$, if 
$G$ intersects $f^{\sharp}$ in a subgroup of finite index, then $G$ contains 
$f^{\sharp}$. 
So, let $G$ be a definable additive subgroup of $L$ such that $[f^{\sharp}:G\cap 
f^{\sharp}]$ is finite. 

Assume that there exists $n\ge 1$ such that 
$\left[f^n(L):G\cap f^n(L)\right]$ is finite. For such $n$, 
because $f^n(L)$ is connected (see Corollary ~\ref{C:fksepconnected}), we 
conclude that $f^n(L)=G\cap f^n(L)$. So, $f^n(L)\subset G$. Then, by the 
definition of $f^{\sharp}$, we get that $f^{\sharp}\subset G$.

Suppose that for all $n\ge 1$, $\left[f^n(L):G\cap f^n(L)\right]$ is infinite. 
By compactness and the fact that $L$ is $\aleph_1$-saturated, we conclude that 
also 
$\left[f^{\sharp}:G\cap f^{\sharp}\right]$ is infinite, which contradicts our 
assumption. 
\end{proof}

\begin{corollary}
\label{L:f,g}
Let $f,g\in K\{\tau_0\}\tau_0\setminus\{0\}$. If $g^{\sharp}\subset f^{\sharp}$, 
then 
$f^{\sharp}=g^{\sharp}$.
\end{corollary}

\begin{proof}
By Theorem \ref{T:fsharpminimal} and our hypothesis, $g^{\sharp}$ is an infinite 
subgroup of $f^{\sharp}$. Thus for every $n\ge 1$, $g^n(L)\cap 
f^{\sharp}$ is a definable infinite subgroup of $f^{\sharp}$. By Theorem 
~\ref{T:fsharpminimal} and Lemma 
~\ref{L:fsharp connected}, $f^{\sharp}\subset 
g^n(L)$. Because this last inclusion holds for all $n\ge 1$, we conclude 
that $f^{\sharp}\subset g^{\sharp}$. Thus $f^{\sharp}=g^{\sharp}$.
\end{proof}

In \cite{Blo} (see Proposition $3.1$ and the \emph{Remark} after the proof of 
Lemma $3.8$) the following result is proved. 
\begin{proposition}
\label{P:lambdapolynomial}
(i) The Frobenius $\tau_0$ and the $\lambda$-functions of level $1$ generate 
$\End(L,+)$ over $K^{\sep}$. Each such element of $\End(L,+)$ will be called an 
(additive) $\lambda$-polynomial. (Because we will only deal with additive 
$\lambda$-polynomials, we will call them simply $\lambda$-polynomials.)

(ii) For every $\psi\in\End(L,+)$, there exists $n\ge 1$ such that for all $g\in 
K^{\sep}\{\tau_0\}\tau_0^n$, $\psi g\in K^{\sep}\{\tau_0\}$.

(iii) Let $G$ be an infinitely definable subgroup of $(L,+)$. Then for each 
$f\in\End(G)$, $f$ extends to an element of $\End(L,+)$.
\end{proposition}

\section{Quasi-endomorphisms of minimal groups associated to Drinfeld modules}
Let $q$ be a 
power of $p$ and let $\tau$ be the power of the Frobenius for which 
$\tau(x)=x^q$, for every $x$. Let $K$ be a finitely generated field extension of 
$\mathbb{F}_q$ of positive transcendence degree. We let $K\{\tau\}$ be the ring 
of all polynomials 
in $\tau$ with coefficients from $K$. Let 
$$f=\sum_{i=0}^r a_i\tau^i\in K\{\tau\},$$
with $a_r\ne 0$. We call the \emph{order} of $f$ at $\tau$ and we denote it by 
$\ord_{\tau}f$, the index $i$ of the first nonzero coefficient $a_i$ of $f$. 
Thus, $f$ is inseparable if and only if $\ord_{\tau}f>0$. 

Let $C$ be a non-singular projective curve defined over $\mathbb{F}_q$. Let $A$ 
be the ring of regular functions on $C$ away from a fixed closed point of $C$. 
Then $A$ is a Dedekind domain. Let $i:A\rightarrow K$ be a morphism. We call the 
morphism $\phi:A\rightarrow K\{\tau\}$ a Drinfeld module if for every $a\in A$, 
the coefficient of $\tau^0$ in $\phi_a$ is $i(a)$, and if there exists $a\in A$ 
such that $\phi_a\ne i(a)\tau^0$. Following the definition from \cite{Goss}, we 
call $\phi$ a Drinfeld module of \emph{generic} characteristic if $\ker 
(i)=\{0\}$. If $\ker(i)=\mathfrak{p}\ne \{0\}$, we call $\phi$ a Drinfeld module 
of 
\emph{finite} characteristic $\mathfrak{p}$. If $\phi$ is a Drinfeld module of 
generic 
characteristic, we let $i:A\rightarrow K$ extend to an embedding of 
$\Frac(A)\subset K$.

As in Section $1$, let $L$ be an $\aleph_1$-saturated elementary extension of 
$K^{\sep}$.

\begin{definition}
\label{D:phisharp}
Let $\phi:A\rightarrow K\{\tau\}$ be a Drinfeld module of finite characteristic. 
We define $$\phi^{\sharp}=\phi^{\sharp}(L)=\cap_{a\in 
A\setminus\{0\}}\phi_a(L).$$
\end{definition}

\begin{lemma}
\label{L:real phisharp}
Let $\phi:A\rightarrow K\{\tau\}$ be a Drinfeld module of finite characteristic 
$\mathfrak{p}$. Let $t\in\mathfrak{p}\setminus\{0\}$. Then 
$$\phi^{\sharp}=\cap_{n\ge 
1}\phi_{t^n}(L)=(\phi_t)^{\sharp}.$$ 
\end{lemma}

\begin{proof}
If $a\notin \pp$, then $\phi_a$ is a separable polynomial and 
$\phi_a(L)=L$. Thus 
\begin{equation}
\label{E:ainrho}
\phi^{\sharp}=\cap_{a\in \pp\setminus\{0\}}\phi_a(L).
\end{equation}
Let $a\in\pp\setminus\{0\}$. Because $t\in\pp\setminus\{0\}$, there exist 
$n,m\ge 1$ and there exist $u,v\in A\setminus \pp$ such that $t^nv=a^mu$. Then 
$\phi_u$ and $\phi_v$ are separable and so,
\begin{equation}
\label{E:uniformizer}
\phi_{a^m}(L)=\phi_{a^m}(\phi_u(L))=\phi_{a^mu}(L)=\phi_{t^nv}(L)= 
\phi_{t^n}(\phi_v(L))=\phi_{t^n}(L).
\end{equation}
So, $\phi_{t^n}(L)\subset\phi_a(L)$. Thus, using \eqref{E:ainrho}, 
we conclude that the result of Lemma 
~\ref{L:real phisharp} holds.
\end{proof}

The following result is an immediate consequence of Lemmas \ref{L:real phisharp} 
and \ref{L:fsharp connected} and Theorem \ref{T:fsharpminimal}.
\begin{corollary}
\label{C:phisharp is minimal and connected}
The group $\phi^{\sharp}$ is a $c$-minimal, connected additive group.
\end{corollary}

\begin{lemma}
\label{L:endphisubsetqua}
Let $\phi$ be a Drinfeld module of finite characteristic. Then 
$\End(\phi)\subset\End(\phi^{\sharp})\subset\qua(\phi^{\sharp})$.
\end{lemma}

\begin{proof}
Let $t$ be a uniformizer of the prime ideal of $A$ which is the characteristic 
of $\phi$. The inclusion $\End(\phi^{\sharp})\subset\qua(\phi^{\sharp})$ is 
clear. 
Let now $f\in\End(\phi)$ and $x\in\phi^{\sharp}$. We need to show that 
$f(x)\in\phi^{\sharp}$. Because $x\in\phi^{\sharp}$, for all $n\ge 1$, there 
exists $x_n\in L$ such 
that $x=\phi_{t^n}(x_n)$. Because $f\in\End(\phi)$, 
$f(x)=f(\phi_{t^n}(x_n))=\phi_{t^n}(f(x_n))\in\phi_{t^n}(L)$, for all 
$n\ge 1$. Thus indeed, $f(x)\in\phi^{\sharp}$ (see Lemma ~\ref{L:real 
phisharp}).
\end{proof}

\begin{corollary}
\label{C:extended phisharp}
If $\phi$ is a finite characteristic Drinfeld module, then 
$\phi^{\sharp}=\cap_{f\in\End(\phi^{\sharp})}f(L)$.
\end{corollary}

\begin{proof}
For every nonzero $a\in A$, 
$\phi_a\in\End(\phi)\subset\End(\phi^{\sharp})$. Thus
$$\cap_{f\in\End(\phi^{\sharp})}f(L)\subset\phi^{\sharp}.$$
But by Lemma ~\ref{L:minimal} and Corollary ~\ref{C:phisharp is 
minimal and connected}, all the endomorphisms of $\phi^{\sharp}$ are surjective 
on $\phi^{\sharp}$. So, then indeed
$$\phi^{\sharp}=\cap_{f\in\End(\phi^{\sharp})}f(L).$$
\end{proof}

Using Corollary \ref{C:phisharp is minimal and connected} and Proposition 
\ref{P:lambdapolynomial}, we get the following result.
\begin{corollary}
\label{P:inseplambda}
Let $f\in\End(\phi^{\sharp})$. Then $f$ is a $\lambda$-polynomial. In 
particular, there exists $m\ge 1$ such that for all $h\in 
K^{\sep}\{\tau\}\tau^m$, 
$fh\in K^{\sep}\{\tau\}$.
\end{corollary}

For every $a\in A\setminus\{0\}$, we let $\phi[a]=\{x\in K^{\alg}\mid 
\phi_a(x)=0\}$. Then for $a\in A\setminus\{0\}$, we let 
$\phi[a^{\infty}]=\cup_{n\ge 1}\phi[a^n]$. If $\pp$ is any nontrivial prime 
ideal in $A$, then we define
$$\phi[\pp ']=\{x\in K^{\alg}\mid \text{ there exists }a\notin\pp\text{ such 
that }\phi_a(x)=0\}.$$

We define $\phi^{\sharp}(K^{\sep})=\phi^{\sharp}(L)\cap K^{\sep}$. We claim that 
this definition for $\phi^{\sharp}(K^{\sep})$ is equivalent with 
$\phi^{\sharp}(K^{\sep})=\cap_{a\in A\setminus\{0\}}\phi_a(K^{\sep})$. Indeed, 
if $x\in\phi^{\sharp}(L)\cap K^{\sep}$, then for every $a\in A\setminus\{0\}$, 
there exists $x_a\in L$ such that $x=\phi_a(x_a)$. Because $\phi_a\in 
K^{\sep}\{\tau\}$ and $L$ is an elementary extension of $K^{\sep}$, $x_a\in 
K^{\sep}$. Moreover, a similar proof as in Lemma ~\ref{L:real phisharp}, shows 
that $\phi^{\sharp}=\cap_{n\ge 1}\phi_{t^n}(K^{\sep})$, if $\phi_t$ is 
inseparable.

We will continue to denote by $\phi^{\sharp}$ the group $\phi^{\sharp}(L)$ and 
by $\phi^{\sharp}(K^{\sep})$, its subgroup contained in $K^{\sep}$.
\begin{lemma}
\label{L:phisharp is big}
Let $\phi:A\rightarrow K\{\tau\}$ be a Drinfeld module of finite characteristic 
$\pp$. Then $\phi[\pp']\subset\phi^{\sharp}(K^{\sep})$.
\end{lemma}

\begin{proof}
Let $x\in\phi[\pp']$ and let $a\notin\pp$ such that $\phi_a(x)=0$. Because 
$\phi_a$ is separable, $x\in K^{\sep}$. Let $t$ be a uniformizer for $\pp$. 
Because $x\in\phi_{\tor}\setminus\phi[t^{\infty}]$, 
the sequence $\left(\phi_{t^n}(x)\right)_{n\ge 0}$ is periodic. Thus, there 
exists $N_1\ge 1$ such that for all $n\ge 0$, $x=\phi_{t^{nN_1}}(x)$. By Lemma 
~\ref{L:real phisharp}, we conclude that $x\in\phi^{\sharp}(K^{\sep})$. 
\end{proof}

\begin{theorem}
\label{T:t}
Let $\phi:A\rightarrow K\{\tau\}$ be a Drinfeld module of finite characteristic 
$\pp$. Assume there exists a non-constant $t\in A\setminus\{0\}$ such that 
$\phi[t^{\infty}]\cap K^{\sep}$ is finite. Then 
$\phi^{\sharp}(K^{\sep})=\phi[\pp']$. Moreover, with the above hypothesis on 
$\phi_t$, we have that for every $\psi\in\qua(\phi^{\sharp})$, there exists 
$n\ge 1$ such that $\psi\phi_{t^n}=\phi_{t^n}\psi$ (the identity being seen in 
$\qua(\phi^{\sharp})$).
\end{theorem}

\begin{proof}
Clearly, $t\in\pp\setminus\{0\}$, because for all $a\in A\setminus\pp$, 
$\phi_a$ is separable and so, $\phi[a^{\infty}]\subset K^{\sep}$. By Lemma 
\ref{L:real phisharp}, we know that
\begin{equation}
\label{E:realphisharp}
\phi^{\sharp}=\cap_{n\ge 1}\phi_{t^n}(L)
\end{equation}
and $\phi^{\sharp}(K^{\sep})=\cap_{n\ge 1}\phi_{t^n}(K^{\sep})$.

Because $\phi[t^{\infty}]\cap K^{\sep}$ is finite, let $N_0\ge 1$ satisfy
\begin{equation}
\label{E:N_0}
\phi[t^{\infty}]\cap K^{\sep}\subset \phi[t^{N_0}].
\end{equation}
Thus
\begin{equation}
\label{E:zero phisharp phitinfty}
\phi[t^{\infty}]\cap\phi^{\sharp}=\{0\}.
\end{equation}

\begin{lemma}
\label{C:phitaut}
Under the hypothesis of Theorem \ref{T:t}, $\phi_t\in\Aut(\phi^{\sharp})$.
\end{lemma}

\begin{proof}[Proof of Lemma \ref{C:phitaut}.]
By Lemma \ref{L:endphisubsetqua}, we know that $\phi_t\in\End(\phi^{\sharp})$. 
By the definition of $\phi^{\sharp}$, we know that $\phi_t$ is a surjective 
endomorphism of $\phi^{\sharp}$. By \eqref{E:zero phisharp phitinfty}, we know 
that $\phi_t$ is an injective endomorphism of $\phi^{\sharp}$.
\end{proof}

\begin{lemma}
\label{L:coherence}
Assume $x\in\phi^{\sharp}(K^{\sep})$. We can find a sequence $(x_n)_{n\ge 
0}\subset\phi^{\sharp}(K^{\sep})$ such that $x_0=x$ and for all $n\ge 0$, 
$\phi_t(x_{n+1})=x_n$.
\end{lemma}
\begin{proof}[Proof of Lemma \ref{L:coherence}.]
Initially we know that there exists a sequence $(x_n)_{n\ge 0}\subset K^{\sep}$, 
with $x_0=x$ such that for all $n\ge 1$, $x_0=\phi_{t^n}(x^n)$ but we do not 
know that $(x_n)_{n\ge 1}\subset\phi^{\sharp}$ or that $x_n=\phi_t(x_{n+1})$.

\emph{Claim} For every $x\in\phi^{\sharp}(K^{\sep})$, we can always choose 
arbitrarily long (finite) sequences $$(x_n)_{0\le n\le N}\subset K^{\sep}$$ such 
that $x_0=x$ and for every $0\le n<N$, $x_n=\phi_t(x_{n+1})$. 

\begin{proof}[Proof of Claim.]
Because $x\in\phi^{\sharp}(K^{\sep})$, pick $x_N\in K^{\sep}$ such that 
$x=\phi_{t^N}(x_N)$ and then for $0\le n\le N-1$, define 
$x_n=\phi_{t^{N-n}}(x_N)$.
\end{proof}

\begin{sublemma}
\label{S:S1}
For every $x\in\phi^{\sharp}(K^{\sep})$, there exists 
$x_1\in\phi^{\sharp}(K^{\sep})$ such that $x=\phi_t(x_1)$.
\end{sublemma}

\begin{proof}[Proof of Sublemma \ref{S:S1}.]
By the result obtained in the above \emph{Claim}, we can construct a coherent 
sequence 
$(x_n)_{0\le n\le N_0+1}$ (for $N_0$ as in \eqref{E:N_0}) corresponding to $x$. 
Then $\phi_t(x_1)=x$ and we claim that $x_1\in\phi^{\sharp}(K^{\sep})$. For this 
we need 
to show that for every $n\ge 1$, there exists $z_n\in K^{\sep}$ such that 
$x_1=\phi_{t^n}(z_n)$. We will prove this assertion by induction on $n$.

We already know that the statement holds for $n\le N_0$, because $x_1$ is part 
of the the coherent sequence $(x_n)_{0\le n\le N_0+1}$. So, we only need to show 
the inductive step. Thus, we 
suppose that for some $n\ge N_0$, there exists $z_n\in K^{\sep}$ such that 
$x_1=\phi_{t^n}(z_n)$ and we will prove that there exists $z_{n+1}\in K^{\sep}$ 
such that $x_1=\phi_{t^{n+1}}(z_{n+1})$.

Because $x\in\phi^{\sharp}(K^{\sep})$, there exists $x_{n+2}\in K^{\sep}$ such 
that 
$x=\phi_{t^{n+2}}(x_{n+2})$. But we already know that $x=\phi_t(x_1)$ and by the 
induction hypothesis we also know that $x_1=\phi_{t^n}(z_n)$. Thus 
$x=\phi_{t^{n+1}}(z_n)$ and so, 
\begin{equation}
\label{E:8}
x=\phi_{t^{n+1}}(z_n)=\phi_{t^{n+2}}(x_{n+2}).
\end{equation}
Equation \eqref{E:8} implies that 
\begin{equation}
\label{E:9}
z_n-\phi_t(x_{n+2})\in\phi[t^{n+1}].
\end{equation}
We know that $z_n$ and $x_{n+2}$ are in $K^{\sep}$, while $\phi_t\in K\{\tau\}$ 
and so, 
\begin{equation}
\label{E:91}
\omega=z_n-\phi_t(x_{n+2})\in\phi[t^{n+1}]\cap K^{\sep}.
\end{equation}
But because of \eqref{E:N_0}, it means that actually 
\begin{equation}
\label{E:10}
\omega\in\phi[t^{N_0}]\cap K^{\sep}.
\end{equation}
Thus $\phi_{t^{N_0}}(\omega)=0$ and because $n\ge N_0$ it means that also 
$\phi_{t^n}(\omega)=0$. So, using equation \eqref{E:9} we get that 
$\phi_{t^n}(z_n)=\phi_{t^n}(\phi_t(x_{n+2}))$. But by our induction hypothesis, 
$x_1=\phi_{t^n}(z_n)$ and so we get that $x_1=\phi_{t^{n+1}}(x_{n+2})$. Thus, 
for $z_{n+1}=x_{n+2}$ the inductive step holds.

So, indeed $x_1\in\phi^{\sharp}(K^{\sep})$ and $x=\phi_t(x_1)$
\end{proof}

Using the result of \eqref{S:S1} we are able to conclude the proof of Lemma 
\ref{L:coherence}. For $x\in\phi^{\sharp}(K^{\sep})$ we find as in Sublemma 
\ref{S:S1}, 
$x_1\in\phi^{\sharp}(K^{\sep})$ such that $x=x_0=\phi_t(x_1)$. Then we apply the 
result of 
\eqref{S:S1} to $x_1$ and find $x_2\in\phi^{\sharp}(K^{\sep})$ such that 
$x_1=\phi_t(x_2)$ 
and we continue this process until we construct the desired infinite coherent 
sequence $(x_n)_{n\ge 0}$ for which $x_{n}=\phi_t(x_{n+1})$, for all $n$.
\end{proof}

\begin{proof}[Alternative proof of Lemma \ref{L:coherence}.] Let 
$x\in\phi^{\sharp}(K^{\sep})$. As shown by the \emph{Claim} inside the proof of 
Lemma ~\ref{L:coherence}, there are arbitrarily long (finite) sequences 
$(x_n)_{0\le n\le N}$ such that $x=x_0$ and for every $0\le n\le N-1$, 
$x_n=\phi_t(x_{n+1})$. By compactness, because $L$ is $\aleph_1$-saturated, 
there exists an infinite coherent sequence $(x_n)_{n\ge 0}\subset L$ such that 
$x=x_0$ and for every $n\ge 0$, $x_n=\phi_t(x_{n+1})$. Because $x\in K^{\sep}$ 
and $\phi_t\in K\{\tau\}$, $(x_n)_{n\ge 0}\subset K^{\alg}\cap L=K^{\sep}$ (the 
intersection of the two fields being taken inside a fixed algebraic closure of 
$L$ which contains $K^{\alg}$). 
\end{proof}

The result of Lemma \ref{L:coherence} is instrumental in proving that 
$\phi^{\sharp}(K^{\sep})\subset\phi_{\tor}$. Indeed, take 
$x\in\phi^{\sharp}(K^{\sep})$ and 
construct the associated sequence $(x_n)_{n\ge 0}$ as in \eqref{L:coherence}.

Let $K'=K(x)$. We claim that $x_n\in K'$, for all $n\ge 1$.

Fix $n\ge 1$ and pick any $\sigma\in\Gal(K^{\sep}/K')$. Because $\phi_t\in 
K\{\tau\}\subset K'\{\tau\}$, for every 
$m\ge 1$, $\sigma(x_m)=\sigma(\phi_t(x_{m+1}))=\phi_t(\sigma(x_{m+1}))$. So, for 
every $m\ge 1$, $x_n-\sigma(x_n)=\phi_{t^m}(x_{n+m}-\sigma(x_{n+m}))$. Thus, 
\begin{equation}
\label{E:12}
x_{n}-\sigma(x_{n})\in\phi^{\sharp}.
\end{equation}
But 
$\phi_{t^n}(x_n-\sigma(x_n))=\phi_{t^n}(x_n)-\phi_{t^n}(\sigma(x_n))=\phi_{t^n}(
x_n)-\sigma(\phi_{t^n}(x_n))=x-\sigma(x)=0$, because $x\in K'$. Thus
\begin{equation}
\label{E:11}
x_{n}-\sigma(x_{n})\in\phi[t^n].
\end{equation}
As shown by \eqref{E:zero phisharp phitinfty}, there is no $t$-power torsion of 
$\phi$ in $\phi^{\sharp}$. Equations \eqref{E:11} and \eqref{E:12} yield 
\begin{equation}
\label{E:13}
x_{n}-\sigma(x_{n})=0.
\end{equation}
So, $x_n=\sigma(x_n)$, for all $n\ge 1$ and for all 
$\sigma\in\Gal(K^{\sep}/K')$. 
Thus, $x_n\in K'$, for all $n\ge 1$ as it was claimed. If $x\notin\phi_{\tor}$ 
it 
means that $x_n\notin\phi_{\tor}$, for all $n\ge 1$. This will give a 
contradiction to the structure theorem for $\phi(K')$. 

In \cite{Poo} (for fields of transcendence degree $1$ over 
$\mathbb{F}_p$) and 
in \cite{Wan} (for fields of arbitrary positive transcendence degree) it is 
established 
that a finitely generated field (as $K'$ in our setting) has the following 
$\phi$-module structure: a direct sum of a finite torsion submodule and a free 
module of rank $\aleph_0$. In particular this means that there cannot be an 
infinitely $t$-divisible non-torsion element $x\in L$. So, $x\in\phi_{\tor}$ and 
we conclude that $\phi^{\sharp}(K^{\sep})\subset\phi_{\tor}$. 

By Lemma \ref{L:phisharp is big}, we know that 
$\phi[\pp']\subset\phi^{\sharp}$. We will prove next that under the hypothesis 
from Theorem ~\ref{T:t} (see \eqref{E:N_0}), 
$\phi^{\sharp}(K^{\sep})=\phi[\pp']$.

Suppose that there exists $x\in\phi^{\sharp}(K^{\sep})\setminus\phi[\pp']$. 
Because we 
already proved that $\phi^{\sharp}(K^{\sep})\subset\phi_{\tor}$, 
$x\in\phi_{\tor}$. Then 
there exists $a\in\pp\setminus\{0\}$ such that $\phi_a(x)=0$. Because 
$t\in\pp\setminus\{0\}$, there exist $n,m\ge 1$ and $u,v\in A\setminus\pp$ 
such that $t^nv=a^mu$. Then
$$\phi_{t^nv}(x)=\phi_{a^mu}(x)=\phi_{a^{m-1}u}(\phi_a(x))=0.$$
So, $x\in\phi[t^nv]$. By our assumption, $x\notin\phi[\pp']$ and so, 
$y=\phi_v(x)\ne 0$. Thus 
\begin{equation}
\label{E:nonzero phitinfty}
y\in\phi[t^{n}]\setminus\{0\}.
\end{equation}
By Lemma ~\ref{L:endphisubsetqua}, because $x\in\phi^{\sharp}(K^{\sep})$ and 
$\phi_v\in\End(\phi)$,
\begin{equation}
\label{E:nonzero phisharp}
y=\phi_v(x)\in\phi^{\sharp}(K^{\sep}).
\end{equation}
Equations \eqref{E:nonzero phitinfty} and \eqref{E:nonzero phisharp} provide a 
contradiction to \eqref{E:zero phisharp phitinfty}. So, indeed 
$\phi^{\sharp}(K^{\sep})=\phi[\pp']$.

In order to prove the second part of our Theorem \ref{T:t} regarding the 
quasi-endomorphisms of $\phi^{\sharp}$, we split the proof in two cases.

\emph{Case 1.} The polynomial $\phi_t$ is purely inseparable.

Then $\phi_t=\alpha\tau^r$ for some $\alpha\in K$ and some $r\ge 1$. Let 
$\gamma\in K^{\sep}$ such that $\gamma^{q^r-1}\alpha=1$. 

Let $\phi^{(\gamma)}$ be the Drinfeld module defined by 
$\phi^{(\gamma)}=\gamma^{-1}\phi\gamma$. We call $\phi^{(\gamma)}$ the conjugate 
of $\phi$ by $\gamma$. Then $\phi^{(\gamma)}_t=\tau^r$. Moreover, because for all 
$a\in A$, $\phi^{(\gamma)}=\gamma^{-1}\phi_a\gamma$ and $\gamma\in K^{\sep}$, we 
conclude that
\begin{equation}
\label{E:conjugated phisharp}
\phi^{(\gamma)^{\sharp}}=\gamma^{-1}\phi^{\sharp}
\end{equation}
and 
\begin{equation}
\label{E:moreover}
\qua(\phi^{\sharp})=\gamma\qua\left(\phi^{(\gamma)^{\sharp}}\right)\gamma^{-1}.
\end{equation}

Because $\phi_t^{(\gamma)}=\tau^r$, 
$\phi^{(\gamma)^{\sharp}}=\cap_{n\ge 1}L^{p^n}:=L^{p^{\infty}}$. By \cite{Blo} 
(Proposition $4.10$), 
$\qua\left(L^{p^{\infty}}\right)$ is the unique field of fractions of 
the Ore's ring $\mathbb{F}_p^{\alg}\{\tau_0,\tau_0^{-1}\}$, where $\tau_0$ is 
the usual Frobenius. 
Then clearly, for all $\psi\in\qua\left(\phi^{(\gamma)^{\sharp}}\right)$, there 
exists $n\ge 1$ such that $\phi^{(\gamma)}_{t^n}$ commutes with $\psi$ in 
$\qua(\phi^{(\gamma)^{\sharp}})$. By \eqref{E:moreover}, we conclude that also 
for every $\psi\in\qua(\phi^{\sharp})$, there exists $n\ge 1$ such that 
$\psi\phi_{t^n}=\phi_{t^n}\psi$.

\emph{Case 2.} The polynomial $\phi_t$ is not purely inseparable, i.e. 
$\phi[t]\ne\{0\}$.

\begin{lemma}
\label{L:qua=end}
For every $\psi\in\qua(\phi^{\sharp})$ there exists $a\in A\setminus\{0\}$ and 
$n\ge 1$ such that 
$\phi_a\psi\phi_{t^n}\in\End(\phi^{\sharp})\cap K^{\sep}\{\tau\}$.
\end{lemma}
\begin{proof}
Let $\psi\in\qua(\phi^{\sharp})$ and let $S=\{x\in\phi^{\sharp}|(0,x)\in\psi\}$. 
Thus, $S$ is a finite, $K^{\sep}$-definable subgroup of $\phi^{\sharp}$. Because 
$L$ is an elementary extension of $K^{\sep}$, $S\subset K^{\sep}$. Thus 
$S\subset\phi^{\sharp}(K^{\sep})\subset\phi_{\tor}$. Hence there exists $a\in 
A\setminus\{0\}$ such that $S\subset\phi[a]$. By 
Lemma ~\ref{L:endphisubsetqua}, $\phi_a\psi\in\qua(\phi^{\sharp})$ and its 
cokernel is 
trivial by our choice for $a$. Thus, $\phi_a\psi$ is actually an endomorphism of 
$\phi^{\sharp}$. 
Also, according to Proposition ~\ref{P:lambdapolynomial}, the endomorphisms of 
$\phi^{\sharp}$ are $\lambda$-polynomials. Thus, by Corollary 
~\ref{P:inseplambda}, because $\phi_t$ is inseparable, there exists $n\ge 1$ 
such that $\phi_a\psi\phi_{t^n}\in\End(\phi^{\sharp})\cap K^{\sep}\{\tau\}$. 
\end{proof}

\begin{proposition}
\label{P:com}
Let $R$ be a domain, i.e. a unital ring with no nontrivial divisors. 

a) Let $y\in R$ be nonzero and suppose that $g\in R$ commutes with $y$ and $xy$ 
for some $x\in R$. Then $g$ also commutes with $x$.

b) Let $y\in R$ be nonzero and suppose that $g\in R$ commutes with $y$ and $yx$ 
for some $x\in R$. Then $g$ also commutes with $x$.
\end{proposition}

\begin{proof}[Proof of Proposition \ref{P:com}.]
It suffices to prove $a)$, because the proof of $b)$ is almost identical; we 
will need to interchange the order of $x$ and $y$ only.

Thus, for the proof of $a)$, we know that 
\begin{equation}
\label{E:commutation}
(gx)y=g(xy)=(xy)g=x(yg)=x(gy)=(xg)y.
\end{equation}
Because $y\in R\setminus \{0\}$ and $R$ is a domain, equation 
\eqref{E:commutation} concludes the proof of Proposition ~\ref{P:com} $a)$.
\end{proof}

We use Proposition ~\ref{P:com} with $R=\qua(\phi^{\sharp})$ because from 
Proposition ~\ref{P:quadivisionring}, we know that $\qua(\phi^{\sharp})$ is a 
division ring. Then by Lemma ~\ref{L:qua=end} and Proposition ~\ref{P:com}, it 
suffices to prove Theorem ~\ref{T:t} for $f\in\End(\phi^{\sharp})\cap 
K^{\sep}\{\tau\}$. 

Let $f\in\End(\phi^{\sharp})\cap K^{\sep}\{\tau\}$. By Corollary 
~\ref{C:phitaut}, $\phi_t^{-1}\in\End(\phi^{\sharp})$ and so, 
$\phi_t^{-1}f\in\End(\phi^{\sharp})$. Hence, $\phi_t^{-1}f$ is a 
$\lambda$-polynomial. By Proposition ~\ref{P:inseplambda}, there exists $m\ge 1$ 
such that for every polynomial $h\in K\{\tau\}\tau^m$,
\begin{equation}
\label{E:15}
\phi_t^{-1}h\in K^{\sep}\{\tau\}.
\end{equation}

Because $\phi_t$ has inseparable degree at least $1$ and $f\in 
K^{\sep}\{\tau\}$, equation \eqref{E:15} 
yields that $g_1=\phi_t^{-1}f\phi_{t^m}\in K^{\sep}\{\tau\}$. Moreover, because 
of Lemma ~\ref{C:phitaut}, $g_1\in\End(\phi^{\sharp})$. This means that the 
equation
\begin{equation}
\label{E:16}
f\phi_{t^m}=\phi_tg_1,
\end{equation}
which initially was true only on $\phi^{\sharp}$ is an identity in 
$K^{\sep}\{\tau\}$. Indeed, $\phi^{\sharp}$ is infinite (see Lemma 
~\ref{L:phisharp is big}) and so, \eqref{E:16} holds for infinitely many points 
of $L$. Thus, because $f\phi_{t^m}$ and $\phi_t g_1$ are polynomials, 
\eqref{E:16} holds identically in $L$. 

Because in equation \eqref{E:16} all the functions are polynomials in $\tau$, we 
can equate the order of $\tau$ in $g_1$. We obtain
\begin{equation}
\label{E:17}
\ord_{\tau}g_1=\ord_{\tau}f+(m-1)\ord_{\tau}\phi_t\ge (m-1)\ord_{\tau}\phi_t\ge 
m-1.
\end{equation}
Thus $\ord_{\tau}(g_1\phi_t)\ge m$ and using \eqref{E:15}, we get that 
$\phi_t^{-1}g_1\phi_t\in\End(\phi^{\sharp})\cap K^{\sep}\{\tau\}$. So, denote by 
$g_2=\phi_t^{-1}g_1\phi_t$. This means that the identity
\begin{equation}
\label{E:171}
\phi_tg_2=g_1\phi_t,
\end{equation} 
which initially was true only on $\phi^{\sharp}$ is actually true everywhere. It 
is the same argument as above when we explained that equation \eqref{E:16} is an 
identity of polynomials from $K^{\sep}\{\tau\}$. 

We equate the order of $\tau$ of the polynomials from \eqref{E:171} and conclude 
that 
\begin{equation}
\label{E:18}
\ord_{\tau}g_2=\ord_{\tau}g_1\ge m-1.
\end{equation}
So, then again $\ord_{\tau}(g_2\phi_t)\ge m$ and we can apply \eqref{E:15} and 
find a polynomial $$g_3\in K^{\sep}\{\tau\}\cap\End(\phi^{\sharp})\text{ such 
that }\phi_tg_3=g_2\phi_t.$$
Once again $\ord_{\tau}g_3=\ord_{\tau}g_2$ and so the above process can continue 
and we construct an infinite sequence $(g_n)_{n\ge 1}\in 
K^{\sep}\{\tau\}\cap\End(\phi^{\sharp})$ such that for every $n\ge 1$,
\begin{equation}
\label{E:20}
\phi_tg_{n+1}=g_n\phi_t.
\end{equation}
Let $g_0=f\phi_{t^{m-1}}$. Then, using \eqref{E:16}, we conclude that equation 
\eqref{E:20} holds also for $n=0$.

An easy induction will show that for every $k\ge 1$ and for all $n\ge 0$,
\begin{equation}
\label{E:21}
\phi_{t^k}g_{n+k}=g_n\phi_{t^k}.
\end{equation}
Indeed, case $k=1$ is equation \eqref{E:20}. So, we suppose that \eqref{E:21} 
holds for some $k\ge 1$ and for all $n\ge 0$ and we will prove it holds for 
$k+1$ and all $n\ge 0$. By equations \eqref{E:20} and 
\eqref{E:21} we have that
$$\phi_{t^{k+1}}g_{n+k+1}=\phi_t(\phi_{t^k}g_{n+1+k})=\phi_tg_{n+1}\phi_{t^k}= 
g_n\phi_t\phi_{t^k}=g_n\phi_{t^{k+1}},$$
which proves the inductive step of our assertion. 

Equation \eqref{E:21} shows that for every $k\ge 1$, $g_{n+k}$ maps $\phi[t^k]$ 
into itself, for every $n\ge 0$. Equation \eqref{E:20} shows that all the 
polynomials $g_n$ have the same degree, call it $d$. Because $\phi_t$ is not 
purely inseparable, we may choose $k_0\ge 1$ such that 
\begin{equation}
\label{E:22}
|\phi[t^{k_0}]|>d.
\end{equation}
Because $\phi[t^{k_0}]$ is a finite set and our sequence of 
polynomials $(g_n)_{n\ge 0}$ is infinite, it means that there exist $n_2> n_1\ge 
0$ such that 
\begin{equation}
\label{E:23}
g_{n_1+k_0}|_{\phi[t^{k_0}]}=g_{n_2+k_0}|_{\phi[t^{k_0}]}.
\end{equation}
By another application of the fact that all $g_n$ are polynomials, equations 
\eqref{E:22} and \eqref{E:23} yield that
\begin{equation}
\label{E:25}
g_{n_1+k_0}=g_{n_2+k_0}.
\end{equation}
But then, using \eqref{E:21} (with $k=n_2-n_1$ and $n=n_1+k_0$) we conclude that 
\begin{equation}
\label{E:24}
\phi_{t^{n_2-n_1}}g_{n_2+k_0}=g_{n_1+k_0}\phi_{t^{n_2-n_1}}.
\end{equation}
If we denote by $g$ the polynomial represented by both $g_{n_2+k_0}$ and 
$g_{n_1+k_0}$ (according to \eqref{E:25}), equation \eqref{E:24} shows that $g$ 
commutes with $\phi_{t^{n_2-n_1}}$. We let $n_0=n_2-n_1\ge 1$ and so,
\begin{equation}
\label{E:26}
g\phi_{t^{n_0}}=\phi_{t^{n_0}}g.
\end{equation}
The definition of $g=g_{n_1+k_0}$ and equation \eqref{E:21} (with $k=n_1+k_0$ 
and $n=0$) give
\begin{equation}
\label{E:27}
\phi_{t^{n_1+k_0}}g=g_0\phi_{t^{n_1+k_0}}.
\end{equation}
Equation \eqref{E:26} shows that $\phi_{t^{n_0}}$ commutes with 
$\phi_{t^{n_1+k_0}}g$. Thus, by equation \eqref{E:27}, $\phi_{t^{n_0}}$ commutes 
also with $g_0\phi_{t^{n_1+k_0}}$. We apply now Proposition ~\ref{P:com} $a)$ to 
conclude that $\phi_{t^{n_0}}$ commutes with $g_0$. Because 
$g_0=f\phi_{t^{m-1}}$, another 
application of the above mentioned proposition gives us 
$$\phi_{t^{n_0}}f=f\phi_{t^{n_0}}.$$
\end{proof}

\begin{theorem}
\label{T:te}
Let $\phi$ be a Drinfeld module of finite characteristic $\pp$. Assume that 
there exists $f\in\Aut(\phi^{\sharp})\cap K^{\sep}\{\tau\}\tau$. Then 
$\phi^{\sharp}\subset\phi_{\tor}$ and for all $\psi\in\qua(\phi^{\sharp})$, 
there exists $n\ge 1$ such that $\psi f^n=f^n\psi$ (the identity being seen in 
$\qua(\phi^{\sharp})$).
\end{theorem}

\begin{proof}
Construct another Drinfeld module $\phi':\mathbb{F}_q[t]\rightarrow 
K^{\sep}\{\tau\}$ by $\phi'_t=f$. By Lemma ~\ref{L:real phisharp}, 
$\phi'^{\sharp}=f^{\sharp}$. Using Corollary ~\ref{C:extended phisharp} and 
$f\in\End(\phi^{\sharp})$, we 
get that 
\begin{equation}
\label{E:31}
\phi^{\sharp}\subset\phi'^{\sharp}.
\end{equation}
Because both $\phi^{\sharp}$ and $\phi'^{\sharp}$ are connected, $c$-minimal 
groups (see Corollary
~\ref{C:phisharp is minimal and connected}), we conclude that they are equal.

Because $\phi'_t\in\Aut(\phi^{\sharp})=\Aut(\phi'^{\sharp})$, 
$\phi'[t^{\infty}]\cap K^{\sep}$ is finite and so, we are in the hypothesis of 
Theorem ~\ref{T:t} with $\phi'$ and $t$. Thus, we conclude that
\begin{equation}
\label{E:41}
\phi'^{\sharp}(K^{\sep})=\phi'[(t)'],
\end{equation}
where by $\phi'[(t)']$ we denoted the prime-to-$t$-torsion of $\phi'$.

Because for all $a\in A$, $\phi_a\in\End(\phi)\subset 
\End(\phi^{\sharp})=\End(\phi'^{\sharp})$, by Theorem \ref{T:t}, there exists 
$n_a\ge 1$ such that $\phi_af^{n_a}=f^{n_a}\phi_a$. Because $A$ is finitely 
generated as a $\mathbb{F}_q$-algebra, we can find $n_0\ge 1$ such that for all 
$a\in A$, $\phi_af^{n_0}=f^{n_0}\phi_a$, i.e. $f^{n_0}\in\End(\phi)$.

\emph{Claim} Let $c(t)\in\mathbb{F}_q[t]\setminus\{0\}$ and let $m\ge 1$. Then 
there exists $d(t)\in\mathbb{F}_q[t^m]\setminus\{0\}$ such that $c(t)$ divides 
$d(t)$.

\begin{proof}[Proof of \emph{Claim}.]
Because $\mathbb{F}_q[t]/(c(t))$ is finite and because $\mathbb{F}_q[t^m]$ is 
infinite, there exist $d_1(t)\ne d_2(t)$, both polynomials in 
$\mathbb{F}_q[t^m]$ such that $c(t)$ divides $d(t)=d_1(t)-d_2(t)$. 
\end{proof}

Let $x\in\phi'_{\tor}$ and let $c(t)\in\mathbb{F}_q[t]\setminus\{0\}$ such that 
$\phi'_{c(t)}(x)=0$. By the above \emph{Claim}, we may assume that 
$c(t)\in\mathbb{F}_q[t^{n_0}]$. Because $\phi'_{t^{n_0}}=f^{n_0}\in\End(\phi)$, 
$\phi'_{c(t)}\in\End(\phi)$. 

Let $a$ be a non-constant element of $A$. Then for all $y\in\phi'[c(t)]$, 
$$\phi'_{c(t)}(\phi_a(y))=\phi_a(\phi'_{c(t)}(y))=0.$$
Thus $\phi_a(y)\in\phi'[c(t)]$ for all $y\in\phi'[c(t)]$. Similarly, 
$\phi_{a^m}$ maps $\phi'[c(t)]$ into itself for every $m\ge 1$. Because 
$\phi'[c(t)]$ is a finite set and $x\in\phi'[c(t)]$, there exist $m_2>m_1\ge 1$ 
such that 
$\phi_{a^{m_2}}(x)=\phi_{a^{m_1}}(x)$. Thus $x\in\phi[a^{m_2}-a^{m_1}]$ and 
$a^{m_2}-a^{m_1}\ne 0$ ($a$ is not constant). This shows that $x\in\phi_{\tor}$ 
and because $x$ was an arbitrary torsion point of $\phi'$, it means that 
$\phi'_{\tor}\subset\phi_{\tor}$. Actually, because the above argument can be 
used reversely by starting with an arbitrary torsion point $x$ of $\phi$ and 
concluding that $x\in\phi'_{\tor}$, we have $\phi_{\tor}=\phi'_{\tor}$. In any 
case, the inclusion $\phi'_{\tor}\subset\phi_{\tor}$ is sufficient to conclude 
that 
$$\phi^{\sharp}(K^{\sep})=\phi'^{\sharp}(K^{\sep})\subset\phi'_{\tor}\subset\phi
_{\tor}.$$

Also, Theorem \ref{T:t} applied to $\phi'$ and $f=\phi'_t$ shows that for all 
$\psi\in\qua(\phi'^{\sharp})=\qua(\phi^{\sharp})$, there exists $n\ge 1$ such 
that $\psi f^n=f^n\psi$ (the identity being seen in $\qua(\phi^{\sharp})$).
\end{proof}

The following example shows that we cannot strengthen Theorem \ref{T:t} and 
also shows how Theorem \ref{T:te} applies when we do not have the hypothesis 
of \eqref{T:t}.

\begin{example}
\label{E:lambda}
Assume $p>2$. Let $f=t\tau+\tau^3$. Then, for all $\lambda\in\mathbb{F}_{q^2}$,
\begin{equation}
\label{E:33}
f\lambda=\lambda^qf
\end{equation}
where $\lambda$ is seen as the operator $\lambda\tau^0$.

Define $\phi:\mathbb{F}_q[t]\rightarrow\mathbb{F}_q(t)\{\tau\}$ by 
$\phi_t=f(\tau^0+f)$. We let $K=\mathbb{F}_q(t)$ and we claim that 
\begin{equation}
\label{E:32}
\phi[t^{\infty}]\cap K^{\sep}\text{ is infinite.}
\end{equation}
Because for all $n\ge 1$, $\phi_{t^n}=f^n(\tau^0+f)^n$, 
$\Ker((\tau^0+f)^n)\subset\Ker\phi_{t^n}$. Because $\tau^0+f$ is a separable 
polynomial, all the roots of $(\tau^0+f)^n$ are separable. So, indeed, 
\eqref{E:32} holds.

Equation \eqref{E:32} shows that the hypothesis of Theorem ~\ref{T:t} fails for 
$\phi$ and $t$. We will prove that also the conclusion of Theorem ~\ref{T:t} 
regarding the quasi-endomorphisms of $\phi^{\sharp}$ fails, i.e. there exists a 
quasi-endomorphism of $\phi^{\sharp}$ that does not commute with any power of 
$\phi_t$.

Let $\lambda\in\mathbb{F}_{q^2}\setminus\mathbb{F}_q$. Applying 
Lemma ~\ref{L:real phisharp}, we get that $\phi^{\sharp}=(\phi_t)^{\sharp}$. 
Applying Lemma ~\ref{L:f,g} to $\phi_t$ and $f^2$ we conclude that 
$\phi^{\sharp}=(f^2)^{\sharp}$. But
\begin{equation} 
\label{E:34}
f^2\lambda=\lambda f^2 \text{ (apply twice equation \eqref{E:33}).}
\end{equation}
Thus, with the help of Lemma \ref{L:endphisubsetqua} applied to the Drinfeld 
module $\psi:\mathbb{F}_q[t]\rightarrow K\{\tau\}$ given by $\psi_t=f^2$, we get 
that $$\lambda\in\End(\psi^{\sharp})=\End\left((f^2)^{\sharp}\right)= 
\End(\phi^{\sharp}).$$

Suppose that there exists $n\ge 1$ such that 
$\phi_{t^n}\lambda=\lambda\phi_{t^n}$ on $\phi^{\sharp}$. Because 
$\phi^{\sharp}$ is infinite, $\phi_{t^n}\lambda=\lambda\phi_{t^n}$, as 
polynomials. Then also $\phi_{t^{2n}}\lambda=\lambda 
\phi_{t^{2n}}$. But $\phi_{t^{2n}}=f^{2n}(\tau^0+f)^{2n}$ and using \eqref{E:34} 
and Proposition ~\ref{P:com} applied to the domain $K\{\tau\}$, we get
\begin{equation}
\label{E:35}
(\tau^0+f)^{2n}\lambda=\lambda(\tau^0+f)^{2n}.
\end{equation}
We will prove that \eqref{E:35} is impossible. Because of the skew commutation 
of $f$ and $\lambda$ as shown in 
equation \eqref{E:33}, the only way for equation \eqref{E:35} to hold is if in 
the expansion of $(\tau^0+f)^{2n}$, all the nonzero terms are even powers of 
$f$. 
Pick $p^l$ be the largest power of $p$ that is less or equal to $2n$. Then 
$\binom{2n}{p^l}\ne 0$ (in $\mathbb{F}_p$) and its corresponding power of $f$ is 
odd. This shows that indeed, \eqref{E:35} cannot hold when $p>2$.

On the other hand, $f\in\End(\phi)$ and the hypothesis of Theorem \ref{T:te} is 
verified for $\phi$ and $f$. Indeed, $f\in\End(\phi^{\sharp})$ and $\Ker(f)\cap 
K^{\sep}=\{0\}$; thus $f\in\Aut(\phi^{\sharp})$. As we can see from equation 
\eqref{E:34}, also the conclusion of \eqref{T:te} holds with $n=2$.

For the case $p=2$ we can construct a similar example by taking $f=t\tau+\tau^4$ 
and define the Drinfeld module 
$\phi:\mathbb{F}_q[t]\rightarrow\mathbb{F}_q(t)\{\tau\}$ by 
$\phi_t=f(\tau^0+f)$. In 
this case, $\lambda\in\mathbb{F}_{q^3}\setminus\mathbb{F}_q$ will play the role 
of the endomorphism of $\phi^{\sharp}$ that commutes with a power of an 
endomorphism of $\phi$, i.e. it commutes with $f^3$, but it does not commute 
with any power of $\phi_t$.
\end{example}

\section{Mordell-Lang conjecture for Drinfeld modules}
In this section, we will use the notation $\overline{X}$ for the Zariski closure 
of the variety $X$. Let $K$ be a finitely generated field of positive 
transcendence degree over $\mathbb{F}_p$.

In \cite{Den}, Laurent Denis formulated an analogue of the 
Mordell-Lang Conjecture in the context of Drinfeld modules. Even though the 
formulation from \cite{Den} is for Drinfeld modules of generic characteristic, 
we can ask the same question for Drinfeld modules of finite characteristic. 
Thus, our Statement \ref{C:Con} will cover both cases. Before stating 
\eqref{C:Con} we need a definition.
\begin{definition}
Let $\phi:A\rightarrow K\{\tau\}$ be a Drinfeld module. For $g\ge 0$ we 
consider $\phi$ acting diagonally on $\mathbb{G}_a^g$. An algebraic 
$\phi$-submodule of $\mathbb{G}_a^g$ is a connected algebraic subgroup of 
$\mathbb{G}_a^g$ which is stable under the action of $\phi$.
\end{definition}

\begin{statement}[Mordell-Lang statement for $\phi$]
\label{C:Con}
Let $\phi$ be a Drinfeld module. If $\Gamma$ is a 
finitely generated $\phi$-submodule of $\mathbb{G}_a^{g}(K^{\alg})$ for some 
$g\ge 0$ and if $X$ is an algebraic subvariety of $\mathbb{G}_a^{g}$, then there 
are finitely 
many algebraic $\phi$-submodules $B_1,\dots,B_s$ and there are 
finitely many elements $\gamma_1,\dots,\gamma_s$ of $\mathbb{G}_a^{g}(K^{\alg})$ 
such that 
$X(K^{\alg})\cap\Gamma =\cup_{1\le i\le s}(\gamma_i+B_i(K^{\alg})\cap\Gamma)$.
\end{statement}

The first result towards Conjecture \ref{C:Con} was obtained by Thomas Scanlon 
in \cite{Sca}. Before stating the theorem from \cite{Sca}, we need to introduce 
two 
definitions.

\begin{definition}
\label{D:fieldofdef}
For a Drinfeld module $\phi:A\rightarrow K\{\tau\}$, its field of definition is 
the smallest subfield of $K$ containing all the coefficients of $\phi_a$, for 
every $a\in A$.
\end{definition}

\begin{definition}
\label{D:D1}
Let $\phi:A\rightarrow K\{\tau\}$ be a Drinfeld module. The modular 
transcendence degree of $\phi$ is the minimum transcendence degree over 
$\mathbb{F}_p$ of the field of definition for $\phi^{(\gamma)}$, where the 
minimum is taken over all $\gamma\in K^{\alg}\setminus\{0\}$.
\end{definition}
In \cite{DG2} we proved that if there exists a non-constant $t\in A$ such that 
$\phi_t=\sum_{i=0}^ra_i\tau^i$ is monic, then the modular transcendence degree 
of $\phi$ is $\trdeg_{\mathbb{F}_p}\mathbb{F}_p(a_0,\dots,a_{r-1})$.

\begin{theorem}[Thomas Scanlon]
\label{T:TS}
Let $\phi:A\rightarrow K\{\tau\}$ be a Drinfeld module of finite characteristic 
and modular transcendence degree at least $1$. Let $\Gamma$ be a finitely 
generated $\phi$-submodule of $\mathbb{G}_a^g(K^{\alg})$ and $X$ be an algebraic 
subvariety of $\mathbb{G}_a^g$. Then $X(K^{\alg})\cap\Gamma$ is a finite 
union of translates of subgroups of $\Gamma$.
\end{theorem}

Using Theorem \ref{T:t}, we are able to strengthen the conclusion of 
\eqref{T:TS} by showing that one could replace \emph{subgroups} by 
\emph{$\phi$-submodules}.
\begin{theorem}
\label{T:DG}
If $X$ is an algebraic subvariety of $\mathbb{G}_a^g$ and $\phi:A\rightarrow 
K\{\tau\}$ is a Drinfeld module of positive modular transcendence degree for which there exists a non-constant $t\in A$ 
such that
$\phi[t^{\infty}](K^{\sep})$ is finite, then there exists $n_0\ge 1$ such that
for every finitely generated $\phi$-submodule $\Gamma$ of 
$\mathbb{G}_a^g(K^{\alg})$, 
$X(K^{\alg})\cap\Gamma$ is a finite union of translates of 
$\mathbb{F}_q[t^{n_0}]$-submodules of $\Gamma$.
\end{theorem}

\begin{proof}
At the expense of replacing $K$ by a finite extension, we may assume that $X$ is 
defined over $K$ and $\Gamma\subset\mathbb{G}_a^g(K)$. Then 
$X(K^{\alg})\cap\Gamma=X(K)\cap\Gamma=X(K^{\sep})\cap\Gamma$.

In \cite{Sca} it was proved Theorem \ref{T:TS} without the mention of $\phi_t$ 
neither in 
the hypothesis nor in the conclusion. Thus we only need to show how 
we can infer that there is an action by a power of $\phi_t$ on the algebraic 
subgroups whose translates form $X(K^{\alg})\cap\Gamma$. 

Let $H$ be an algebraic subgroup of $\mathbb{G}_a^g$ such that for some 
$\gamma\in\mathbb{G}_a^g(K^{\alg})$, 
$$\gamma+H(K^{\alg})\cap\Gamma\subset X(K^{\alg})\cap\Gamma.$$
At the expense of replacing $K$ by a finite extension we may assume that $\gamma\in K$ and $H$ is defined over $K$. Also we may assume that $H(K)\cap\Gamma$ is dense in $H$ (otherwise we replace $H$ by the Zariski closure of $H(K)\cap\Gamma$).

Let $d=\dim H$. If $d=g$, then $H=\mathbb{G}_a^g$ and $H$ is invaried by the 
action of $\phi_t$.

Suppose from now on that $d<g$. Then, without loss of generality, we may suppose 
that for every $x_1,\dots,x_d\in K^{\alg}$, there are at most finitely many and 
at least one tuple $(x_{d+1},\dots,x_g)\in \left(K^{\alg}\right)^{g-d}$ such 
that 
\begin{equation}
\label{E:correspondence g-d}
(x_1,\dots,x_g)\in X(K^{\alg}).
\end{equation}
Let $\pi$ be the correspondence that associates to each $d$-tuple 
$(x_1,\dots,x_d)\in \left(K^{\alg}\right)^d$ the finitely many $(g-d)$-tuples 
$(x_{d+1},\dots,x_g)\in\left(K^{\alg}\right)^{g-d}$ such that 
\eqref{E:correspondence g-d} holds. Because $H$ is an algebraic additive group, 
$\pi$ is an additive $K$-definable correspondence. 

Let $L$ be an $\aleph_1$-saturated elementary extension of $K^{\sep}$. Then 
$H(L)\cap\phi^{\sharp}(L)^g$ is the graph of a quasi-morphism between 
$\phi^{\sharp}(L)^d$ and $\phi^{\sharp}(L)^{g-d}$. By \emph{Lemme} $3.5.3$ of 
\cite{BloT}, 
$$\qh\left(\phi^{\sharp}(L)^d,\phi^{\sharp}(L)^{g-d}\right)\isomto  
M_{d,g-d}\left(\qua(\phi^{\sharp}(L))\right),$$
where by $M_{d,g-d}\left(\qua(\phi^{\sharp}(L))\right)$ we denote the ring of $d\times (g-d)$ matrices over the ring $\qua(\phi^{\sharp}(L))$. Then, by Theorem \ref{T:t}, the image of $\pi$ in 
$\qh\left(\phi^{\sharp}(L)^d,\phi^{\sharp}(L)^{g-d}\right)$ commutes with a 
power of $\phi_t$. Let $n_0\ge 1$ such that 
$$\phi_{t^{n_0}}\left(H(L)\cap\phi^{\sharp}(L)^g\right)=H(L)\cap\phi^{\sharp}(L)
^g.$$
Because $L$ is $\aleph_1$-saturated, by compactness we conclude there exists 
$m\ge 1$ such that
\begin{equation}
\label{E:invaried}
\phi_{t^{n_0}}\left(H(L)\cap\phi_{t^m}(L)^g\right)\subset H(L).
\end{equation}

We know that $H(L)\cap\Gamma$ is Zariski dense in $H$. Thus there exists $\alpha\in L^g$ such that $H(L)\cap\left(\alpha+\phi_{t^m}(L)^g\right)$ is Zariski dense in $H$. But 
$$H(L)\cap\left(\alpha+\phi_{t^m}(L)^g\right)=\beta+\left(H(L)\cap\phi_{t^m}(L)^g \right)$$
for some $\beta\in\alpha+\phi_{t^m}(L)^g$. Because $H(L)\cap\left(\alpha+\phi_{t^m}(L)^g\right)$ is Zariski dense in $H$, we conclude that $H(L)\cap\phi_{t^m}(L)^g$ is Zariski dense in $H$. Thus the set $H(L)\cap\phi_{t^m}(L)^g$ is Zariski dense in $H$ and it is mapped by $\phi_{t^{n_0}}$ inside $H(L)$. Hence $H$ is invaried by $\phi_{t^{n_0}}$.
\end{proof}

\begin{remark}
\label{R:theorem is sharp}
The result of Theorem \ref{T:DG} is sharp in the sense that we can only get some 
$n_0\ge 1$ such that the conclusion of \eqref{T:DG} holds. For example, let the 
Drinfeld module $\phi:\mathbb{F}_q[t]\rightarrow\mathbb{F}_q(t)\{\tau\}$ be 
defined by $\phi_t=\tau+t\tau^3$ and 
$\lambda\in\mathbb{F}_{q^2}\setminus\mathbb{F}_q$. Let $X\subset\mathbb{G}_a^2$ 
be the curve $y=\lambda x$ and let $\Gamma$ be the cyclic $\phi$-submodule of 
$\mathbb{G}_a^2(\mathbb{F}_{q^2}(t))$ generated by $(1,\lambda)$. As shown in 
Example ~\ref{E:lambda}, $\phi_{t^2}\lambda=\lambda\phi_{t^2}$. Thus for every 
$n\ge 1$, $\left(\phi_{t^{2n}}(1),\phi_{t^{2n}}(\lambda)\right)\in 
X(\mathbb{F}_{q^2}(t))$. So, $X(\mathbb{F}_q(t)^{\alg})\cap\Gamma$ is Zariski 
dense in $X$. But $X$ is not invaried by $\phi_t$; $X$ is invaried by 
$\phi_{t^2}$. Hence in this example, Theorem \ref{T:DG} holds with $n_0=2$.
\end{remark}

\begin{remark}
\label{R:important example}
If we drop the hypothesis on $\phi_t$ from Theorem \ref{T:DG} we may lose the 
conclusion, as it is shown by the following example.

Let $p>2$ and let the Drinfeld module $\phi$ and $\lambda$, $X$ and the 
$\phi$-submodule $\Gamma$ be as in Remark ~\ref{R:theorem is sharp}. Let 
$u=t+t^2$. As shown in Example ~\ref{E:lambda}, 
$\phi[u^{\infty}]\cap\mathbb{F}_p(t)^{\sep}$ is infinite and $X$ is not invaried 
by any power of $\phi_u$. But, as shown in Remark ~\ref{R:theorem is sharp}, 
$X(\mathbb{F}_p(t)^{\alg})\cap\Gamma$ is infinite. 
\end{remark}

The above two remarks \ref{R:theorem is sharp} and \ref{R:important example} 
show that the result of Theorem ~\ref{T:DG} is the most we can hope towards 
Statement ~\ref{C:Con} for Drinfeld modules of finite characteristic.

\begin{definition}
\label{D:D2}
Let $\phi:A\rightarrow K\{\tau\}$ be a Drinfeld module. Let $K_0$ be any 
subfield of $K$. Then the relative modular transcendence degree 
of $\phi$ over $K_0$ is the minimum transcendence degree over $K_0$ of the 
compositum field of $K_0$ and the field of definition of $\phi^{(\gamma)}$, 
minimum being taken over all $\gamma\in K^{\alg}\setminus\{0\}$.
\end{definition}

\begin{theorem}
\label{T:T3}
Let $\phi :A\rightarrow K\{\tau\}$ be a Drinfeld module of generic 
characteristic and of relative modular transcendence degree at least $1$ over 
$F=\Frac (A)$. Let $g\ge 0$ and $X$ be an algebraic subvariety of 
$\mathbb{G}_a^g$. Assume that $X$ does not contain a translate of a nontrivial 
connected algebraic subgroup of $\mathbb{G}_a^g$. Then for every finitely 
generated $\phi$-submodule $\Gamma$ of 
$\mathbb{G}_a^g(K^{\alg})$, we have that $ X(K^{\alg})\cap\Gamma$ is finite.
\end{theorem}

\begin{proof}
We let $F^{\alg}$ be the algebraic closure of $F$ inside $K^{\alg}$. For any two 
subextensions of $K^{\alg}$, their compositum is taken inside $K^{\alg}$.

In the beginning we will prove several reduction steps.

\emph{Step 1.}  It suffices to prove Theorem \ref{T:T3} for $\Gamma$ of the form 
$\Gamma_0^g$ where $\Gamma_0$ is a finitely generated $\phi$-submodule of 
$\mathbb{G}_a(K^{\alg})$. Indeed, if we let $\Gamma_0$ be 
the finitely generated $\phi$-submodule of $K^{\alg}$ generated by all the $g$ 
coordinate projections of $\Gamma$ then clearly $\Gamma\subset\Gamma_0^g$. So, 
we suppose that $\Gamma$ has the form $\Gamma_0^g$. To simplify the notation we 
work with a finitely generated $\phi$-submodule $\Gamma$ of 
$\mathbb{G}_a(K^{\alg})$ and prove that 
$X(K^{\alg})\cap\Gamma^g$ is finite.

\emph{Step 2.} Let $t$ be a non-constant element of $A$. Let $\gamma\in 
K^{\alg}$ such that for the Drinfeld module 
$\phi^{(\gamma)}=\gamma^{-1}\phi\gamma$, $\phi^{(\gamma)}_t\text{ is monic.}$
We let $\gamma^{-1}X$ be the variety whose vanishing ideal is composed of 
functions of the form $f\circ\gamma$, where $f$ is in the vanishing ideal of 
$X$ and $\gamma$ is interpreted as the multiplication-by-$\gamma$-map on each 
component of $\mathbb{G}_a^g$. The conclusion of Theorem \ref{T:T3} is 
equivalent with showing that 
$$(\gamma^{-1}X)(K^{\alg})\cap (\gamma^{-1}\Gamma)^g\text{ is finite.}$$
The variety $\gamma^{-1}X$ has the same property as $X$: it does not contain a 
translate of a non-trivial connected algebraic subgroup of $\mathbb{G}_a^g$. The 
group $\gamma^{-1}\Gamma$ is a finitely generated $\phi^{(\gamma)}$-submodule. 
So, it suffices to prove Theorem ~\ref{T:T3} under the extra hypothesis that 
$\phi_t$ is monic. From now on, let 
$$\phi_t=\tau^r+a_{r-1}\tau^{r-1}+\dots+a_0\tau^0.$$
 
\emph{Step 3.} Because $\phi$ is defined over $K$ and $\Gamma$ is a finitely 
generated $\phi$-submodule of $K^{\alg}$ there exists a finite extension 
$L$ of $K$ that contains all the elements of $\Gamma$. We replace $K$ by 
$L$ and then the conclusion of Theorem \ref{T:T3} reads $\vert 
X(K)\cap\Gamma^g\vert <\aleph_0$.

\emph{Step 4.} We define the \emph{division hull} of $\Gamma$, by
$$\overline{\Gamma}=\{\gamma\in K^{\alg}\mid\text{there exists }a\in 
A\setminus\{0\}\text{ such that }\phi_a(\gamma)\in\Gamma\}.$$
In \cite{DG2} we proved the following result (Theorem $5.7$).
\begin{theorem}
\label{T:T4}
Let $F$ be a countable field of characteristic $p$ and let $K$ be a finitely 
generated field over $F$. Let $\phi:A\rightarrow K\{\tau\}$ be a Drinfeld module 
of positive relative modular transcendence degree over $F$. Then for every 
finite extension $L$ of $K$, $\phi(L)$ is a direct sum of a finite torsion 
submodule and a free submodule of rank $\aleph_0$.
\end{theorem}
Using the result of Theorem \ref{T:T4} for $F^{\alg}$, which is countable, and 
for $F^{\alg}K$, which is finitely generated over $F^{\alg}$, and for $\phi$, 
which has positive relative modular transcendence degree over $F^{\alg}$, we 
conclude that $\phi(F^{\alg}K)$ is the direct sum of a finite torsion submodule 
and a free module of rank $\aleph_0$. Thus, because $\overline{\Gamma}$ has 
finite rank, $\overline{\Gamma}\cap F^{\alg}K$ is 
finitely generated. At the expense of replacing $K$ by a finite extension of the 
form $F'K$, where $F'$ is a finite extension of $F$, we 
may assume that $\overline{\Gamma}\cap F^{\alg}K\subset K$. 

\emph{Step 5.} We may replace $\Gamma$ by $\overline{\Gamma}\cap K$, which is 
also a finitely generated $\phi$-submodule that contains $\Gamma$. 

Now, assuming all the reductions from \emph{Step 1-5}, let $x_0\in 
F\setminus\mathbb{F}_p^{\alg}$. 
Because $\phi$ is defined over $K$, but not over $F$, $\trdeg_FK>0$. Thus, let 
$n=\trdeg_FK\ge 1$ and $\{x_1,\dots,x_n\}\subset K$ be a transcendence basis for 
$K/F$.

Let $C$ be the normalization of $\mathbb{P}_{\mathbb{F}_q}^1$ in $F$. Let $V$ be 
the normalization of $\mathbb{P}_{\mathbb{F}_q}^{n+1}$ in $K$. Both $C$ and $V$ 
are projective, normal varieties defined over a finite field. Let 
$\pi:V\rightarrow C$ be the rational map induced by the inclusion of $F$ in $K$. 
By blowing up $V$ we may assume that $\pi$ is a morphism.

At the expense of replacing $F$ by a finite extension $F'$ and replacing $K$ by 
$F'K$, we may assume that the generic fiber of $\pi$ is geometrically 
irreducible.

The irreducible divisors $\bet$ of $V$ are of two types:

(i) \emph{vertical}, in which case $\pi(\bet)=\pp$ is a closed point of $C$.

(ii) \emph{horizontal}, in which case $\pi|_{\bet}:\bet\rightarrow C$ is a dominant map.

For each irreducible divisor $\bet$ of $V$, we let $K_{\bet}$ be the residue 
field of $K$ at $\bet$. For any element $x$ in the valuation ring at $\bet$, 
we let $x_{\bet}$ be the reduction of $x$ at $\bet$. Also, we denote by 
$r_{\bet}$ the reduction map at $\bet$. If all the elements of $\Gamma$ are 
integral at $\bet$, we let
$$\Gamma_{\bet}=\{x_{\bet}\mid x\in \Gamma\}.$$

We say that $\phi$ has \emph{good reduction} at $\bet$ if for all $a\in 
A\setminus\{0\}$, all the coefficients of $\phi_a$ are integral at $\bet$ and 
the leading coefficient of $\phi_a$ is a unit in the valuation ring at $\bet$. 
If $\phi$ has good reduction at $\bet$, then we denote by $\phi^{\bet}$ the 
corresponding reduction. 

Let $S$ be the set of horizontal divisors of $V$ that are the irreducible 
components of 
the poles of the coefficients $a_i$ of $\phi_t$. According to Lemma $4.6$ of 
\cite{DG2}, the set $S$ is the set of horizontal irreducible divisors of $V$ 
that are places of bad reduction for $\phi$.

At the expense of replacing $F$ by a finite extension $F'$ and replacing $K$ by 
$F'K$, we may assume that for each $\gamma\in S$, the generic fiber of $\gamma$ 
is geometrically irreducible. So, from now on we work under the additional two 
assumptions:
\begin{equation}
\label{E:(1)}
\pi:V\rightarrow C\text{ is a morphism whose generic fiber is geometrically 
irreducible.}
\end{equation}
\begin{equation}
\label{E:(2)}
\text{for each $\gamma\in S$, the generic fiber of $\gamma$ is geometrically 
irreducible.}
\end{equation}

The following results are standard (see \cite{EGA}).
\begin{lemma}
\label{L:above}
Because the generic fiber of $\pi:V\rightarrow C$ is geometrically irreducible, 
for all but finitely many closed points $\pp\in C$, $\pi^{-1}(\pp)$ is 
geometrically irreducible.
\end{lemma}

\begin{lemma}
\label{L:gammafiber}
Let $\gamma\in S$. Because the generic fiber of $\gamma$ is geometrically 
irreducible, for all but finitely many closed points $\pp\in C$, 
$\gamma\cap\pi^{-1}(\pp)$ is geometrically irreducible.
\end{lemma}

\begin{lemma}
\label{L:key lemma}
Let $T$ be the set of vertical irreducible divisors $\bet$ of $V$ which satisfy 
the following properties:

a) $\phi$ has good reduction at $\bet$.

b) $\phi^{\bet}$ is a finite characteristic Dinfeld module of positive modular 
transcendence degree.

c) the projective variety $\bet$ is regular in codimension $1$. Moreover the 
irreducible divisors of $\bet$ are the irreducible components of divisors of the 
form 
$\gamma_{\bet}=\gamma\cap\bet$, where $\gamma$ is an effective horizontal 
divisor of $V$.

d) for each $\gamma\in S$, $\gamma_{\bet}$ is geometrically irreducible.

e) for all $x\in\Gamma$, $x$ is integral at $\bet$.

Then the set $T$ is cofinite in the set of all vertical irreducible divisors of 
$V$. 
\end{lemma}

\begin{proof}[Proof of Lemma \ref{L:key lemma}.]
We will show that each of the conditions $a)$-$e)$ is verified by all but 
finitely many vertical irreducible divisors of $V$.

$a)$ There are finitely many irreducible divisors of $V$ that are places of bad 
reduction for $\phi$. So, in particular, there are finitely many irreducible 
vertical divisors of $V$ that do not satisfy $a)$.

$b)$ By the definition of reduction at $\bet$ (which is a place sitting above a 
prime divisor of $A$), $\phi^{\bet}$ is a finite characteristic Drinfeld 
module.

Because $\phi$ has positive relative modular transcendence degree over $F$, 
there exists $a\in A$ and a coefficient $c$ of $\phi_a$ such that $c\notin 
F^{\alg}$. We view $c$ as a rational map from the generic fiber $V_{\eta}$ of 
$V$ to $\mathbb{P}_F^1$. We spread out $c$ to a dominant rational map 
$\tilde{c}:V\rightarrow \mathbb{P}_C^1$, whose generic fiber is $c$. For all but 
finitely many closed points $\pp\in C$, $\tilde{c_{\pp}}$ is not constant. 
According to the result of Lemma ~\ref{L:above}, for all but finitely many 
$\pp$, $\pi^{-1}(\pp)=\bet$ is geometrically irreducible. For such $\bet$, 
we identify $\tilde{c_{\pp}}$ with the reduction of $c$ at the place $\bet$, 
denoted $c_{\bet}$. Thus for all but finitely many irreducible vertical 
divisors $\bet$, $c_{\bet}\notin\mathbb{F}_p^{\alg}$. So, for these divisors 
$\bet$, $\phi^{\bet}$ has positive modular transcendence degree (remember that 
$\phi^{\bet}_t$ is still monic because $\phi_t$ is monic).

$c)$ By construction, all the irreducible divisors of $V$ are projective 
varieties. 

Because $V$ is normal, $V$ is regular in codimension $1$, i.e. there exists a 
closed subset $L\subset V$ of codimension at least $2$ such that $V\setminus L$ 
is regular. Because $\dim V=n+1$, $\dim L\le n-1$. Let $L_1,\dots, L_s$ be the 
irreducible components of $L$. By Lemma ~\ref{L:above}, there is an open subset 
$U\subset C$ such that for all $\pp\in U$, $\pi^{-1}(\pp)$ is geometrically 
irreducible. Thus the set $T_U$ of irreducible divisors $\bet$ of $V$ sitting 
above some $\pp\in U$, is cofinite in the set of all irreducible vertical 
divisors of $V$. Moreover, the divisors from $T_U$ are mutually disjoint. Thus 
there are at most $s$ divisors $\bet$ from $T_U$ that contain at least one of 
$L_1,\dots ,L_s$. We shrink $T_U$ so that we exclude this finite set of divisors 
$\bet$ from $T_U$. Then for all $\bet\in T_U$, $\bet$ is regular in 
codimension $1$. Also, $T_U$ is cofinite in the set of all vertical irreducible 
divisors of $V$.

Let $\bet\in T_U$. Let $\nu$ be an irreducible divisor of $\bet$. By Hilbert's 
Nullstellensatz, we can find a rational function $f$ on $V$ that vanishes on 
$\nu$ but not on $\bet$. Let $Z$ be the divisor of the zeros of $f$ on $V$. 

We first show that $Z$ is a horizontal divisor of $V$. Suppose that $Z$ is 
vertical. Then, because $Z$ intersects non-trivially $\bet$ and $\bet$ is the 
only irreducible vertical divisor of $V$ sitting above $\pp=\pi(\bet)$ 
(remember that the divisors from $T_U$ satisfy that $\pi^{-1}(\pp)$ is 
geometrically irreducible), $\bet$ is an irreducible component of $Z$. But this 
contradicts our assumption that $f$ does not vanish identically on $\bet$.

So, $Z$ is an effective horizontal divisor of $V$ which intersects $\bet$ in the 
divisor $Z_{\bet}$ and $\nu$ is an irreducible component of $Z_{\bet}$. Because 
$\nu$ was an arbitrary irreducible divisor of $\bet$, we conclude that 
condition $c)$ of Lemma ~\ref{L:key lemma} is satisfied for all but finitely 
many vertical irreducible divisors of $V$.

$d)$ This is proved by Lemma \ref{L:gammafiber}.

$e)$ Because $\Gamma$ is finitely generated as a $\phi$-module and $\phi$ has 
good reduction at all but finitely many irreducible divisors, the elements of 
$\Gamma$ are integral at all but finitely many irreducible divisors of $V$.
\end{proof}

\begin{lemma}
\label{L:S nonempty}
The set $S$ is nonempty.
\end{lemma}

\begin{proof}[Proof of Lemma \ref{L:S nonempty}.]
Assume $S$ is empty, i.e. for every $a\in A$ and for each coefficient $c$ of 
$\phi_a$, the pole of $c$ is vertical. Because $\phi$ is not defined over 
$F^{\alg}$, there exists $a\in A$ and a coefficient $c$ of $\phi_a$ such that 
$c\notin F^{\alg}$. But then the pole of $c$ is not vertical, which gives a 
contradiction to our assumption that $S$ is empty.
\end{proof}

\begin{lemma}
\label{L:Claim 1}
For all but finitely many $\bet\in T$, the reduction map $r_{\bet}$ is 
injective on 
$\Gamma_{\tor}$.
\end{lemma}
\begin{proof}[Proof of Lemma \ref{L:Claim 1}.]  
Because $\Gamma_{\tor}$ is finite ($\Gamma$ is finitely generated), only 
finitely 
many $\bet$ from $T$ appear as irreducible components of the divisor of zeros 
for some torsion element of $\Gamma$. 
\end{proof}

\begin{lemma}
\label{L:Claim 2}  
There exists a non-constant $a\in A$ such that for all $\bet\in T$, 
$\Gamma_{\bet}\cap\phi^{\bet}[a]=\{0\}$.
\end{lemma}
\begin{proof}[Proof of Lemma \ref{L:Claim 2}.] 
Let $\bet\in T$. Conditions $c)$ and $d)$ of Lemma \ref{L:key lemma} show that 
there are $|S|$ irreducible 
divisors of $\bet$ that are places of bad reduction for $\phi^{\bet}$; they 
are of the form $\gamma_{\bet}$ for $\gamma\in S$. Because $|S|\ge 1$ (as shown 
by Lemma ~\ref{L:S nonempty}), we use 
the result of Corollary $4.18$ $b)$ from \cite{DG2}, to conclude that for all 
$x\in\phi_{\tor}^{\bet}(K_{\bet})$, there exists a polynomial 
$b(t)\in\mathbb{F}_q[t]$ of degree at most $\frac{r^3+r^2+2r}{2}|S|$ such that 
$\phi^{\bet}_{b(t)}(x)=0$. Because $\Gamma_{\bet}\subset K_{\bet}$, Lemma 
~\ref{L:Claim 2} holds with $a\in\mathbb{F}_q[t]$ being any irreducible 
polynomial of degree greater than 
$\frac{r^3+r^2+2r}{2}|S|$.
\end{proof}

\begin{lemma}
\label{L:Claim 3}  
Let $a$ be a non-constant element of $A$. For almost all $\bet\in T$, 
$r_{\bet}:\Gamma/\phi_a(\Gamma)\rightarrow\Gamma_{\bet}/\phi^{\bet}_a(\Gamma_
{\bet})$ is injective.
\end{lemma}

\begin{proof}[Proof of Lemma \ref{L:Claim 3}.]
Suppose there are infinitely many irreducible divisors $\bet$ for which 
the map in \eqref{L:Claim 3} is not injective. Because $\Gamma/\phi_a(\Gamma)$ 
is finite we can find $x\in \Gamma\setminus\phi_a(\Gamma)$ such that for 
infinitely many $\bet\in T$,  
$x_{\bet}\in\phi^{\bet}_a(\Gamma_{\bet})$. But there are only finitely many 
$y\in 
K^{\alg}$ such that $\phi_a(y)=x$. Thus, there exists $y\in 
K^{\alg}\setminus\Gamma$, a 
solution to the equation 
\begin{equation}
\label{E:y-x}
\phi_a(y)=x,
\end{equation}
such that for infinitely many $\bet\in T$, if $\bet'$ is a place of $K(y)$ 
sitting above $\bet$, $y_{\bet'}\in\Gamma_{\bet}\subset 
K_{\bet}$.

By \emph{Step 5}, $\Gamma=\overline{\Gamma}\cap F^{\alg}K$. Thus, because $y\in 
\left(K^{\alg}\setminus\Gamma\right)\cap\overline{\Gamma}$, 
\begin{equation}
\label{E:fapt}
y\in K^{\alg}\setminus F^{\alg}K.
\end{equation}
 
Let $W$ be the normalization of $V$ in $K(y)$ and let $\psi:W\rightarrow V$ be 
the corresponding morphism. If the generic fiber of $\pi\circ\psi:W\rightarrow 
C$ is not geometrically irreducible, let $F'$ be a finite extension of 
$F$ such that the generic fiber of $W'=W\times_{\mathbb{F}_q}F'$ is 
geometrically irreducible. Denote by $K_1=F'K(y)$. Also, denote by 
$V'=V\times_{\mathbb{F}_q}F'$. Then there 
exists a rational map $f:W'\rightarrow V'$ induced by the natural inclusion of 
$F'K=K_2$ in $K_1$. Let $d=[K_1:K_2]>1$ (remember that $y\notin F^{\alg}K$, as 
shown by \eqref{E:fapt}). 

Let $C'$ be the normalization of $C$ in $F'$. There are two morphisms 
$\pi_1:W'\rightarrow C'$ and $\pi_2:V'\rightarrow C'$ induced by the morphism 
$\pi$. Thus $\pi_1=\pi_2\circ f$. Because the generic fibers of both $\pi_1$ and 
$\pi_2$ are geometrically irreducible, for all but finitely many closed points 
$\pp'\in C'$, $\pi_1'^{-1}(\pp')$ and $\pi_2'^{-1}(\pp')$ are geometrically 
irreducible. Also, for all but finitely many closed points $\pp'\in C'$, the 
rational map $f_{\pp'}:W'_{\pp'}\rightarrow V'_{\pp'}$ has 
degree $d$. Thus, for all but finitely many closed points $\pp'\in C'$, if 
$\bet_1$ and $\bet_2$ are the unique vertical divisors of $W'$ and $V'$, 
respectively, sitting above $\pp'$,
\begin{equation}
\label{E:degree}
[K_{1_{\bet_1}}:K_{2_{\bet_2}}]=d>1.
\end{equation}

Assume $\bet\in T$ has the property that $y_{\bet_1}\in K_{\bet}$. Then 
$K_{2_{\bet_2}}=K_{1_{\bet_1}}$. So, if $r_{\bet}$ is not injective on 
$\Gamma/\phi_a(\Gamma)$, then equation \eqref{E:degree} would be false for the 
divisors $\bet_1$ and $\bet_2$ that lie over $\bet$. Because \eqref{E:degree} 
holds for all but finitely many corresponding vertical divisors $\bet_1$ and 
$\bet_2$, we conclude that lemma \eqref{L:Claim 3} holds.
\end{proof}

Using Lemmas \ref{L:Claim 1}, \ref{L:Claim 2} and \ref{L:Claim 3} we prove the 
following key result.
\begin{lemma}
\label{L:L3.1}
For all but finitely many $\bet\in T$, the reduction 
$\Gamma\rightarrow\Gamma_{\bet}$ is injective.
\end{lemma}

\begin{proof}[Proof of Lemma \ref{L:L3.1}.] 
Shrink $T$ so that all of the three lemmas \ref{L:Claim 1}, \ref{L:Claim 2} and 
\ref{L:Claim 3} hold for $\bet\in T$. Also, let $a$ be as in Lemma 
~\ref{L:Claim 2}.

If $x\in \Gamma\cap\Ker (r_{\bet})$, then by Lemma \ref{L:Claim 3}, 
$x\in\phi_a(\Gamma )$. This means that there exists $x_1\in \Gamma$ such that 
$\phi_a(x_1)=x$. Reducing at $\bet$, we get $\phi^{\bet}_a(x_{1_{\bet}})=0$ 
which by Lemma ~\ref{L:Claim 2} implies that $x_{1_{\bet}}=0$. But then 
applying again ~\ref{L:Claim 3}, this time to $x_1$, we conclude $x_1\in 
\phi_a(\Gamma)$; 
i.e. there exists $x_2\in \Gamma$ such that $x_1=\phi_a(x_2)$.

So, repeating the above process, we end up with $x\in\cap _{n\ge 
1}\phi_{a^n}(\Gamma)=\Gamma_{\tor}$ (because $\Gamma$ is finitely generated). 
But, by Lemma ~\ref{L:Claim 1}, $\Gamma_{\tor}$ 
injects through the reduction at $\bet$. Thus $x=0$ and so 
the proof of Lemma ~\ref{L:L3.1} ends.
\end{proof}

Now, the property \emph{"$X$ does not contain any translate of a nontrivial 
connected algebraic subgroup of $\mathbb{G}_a^g$"} is a definable property as 
shown in Lemma 11 (page 203) of \cite{Bou}. This means that this 
property is inherited by all but finitely many of the special fibers of $X$. 
Coupling this result with Lemma ~\ref{L:L3.1}, we see that for all but 
finitely many irreducible vertical divisors $\bet$ of $V$, the reduction 
of $X$, called $X_{\bet}$ is also a variety that satisfies the same
hypothesis as $X$ and moreover, $\Gamma$ injects through such 
reduction. This means that 
\begin{equation}
\label{E:1900}
\vert X(K)\cap \Gamma^g \vert\le\vert 
X_{\bet}(K_{\bet})\cap\Gamma_{\bet}^g\vert.
\end{equation}
According to condition $b)$ of Lemma ~\ref{L:key lemma}, for all $\bet\in T$, 
$\phi^{\bet}$ satisfies the hypothesis of Theorem ~\ref{T:TS}. Thus, applying 
Theorem ~\ref{T:TS}, $X_{\bet}\cap\Gamma_{\bet}$ is a finite union of 
translates of cosets of 
subgroups of $\Gamma_{\bet}$. Suppose that one of these subgroups of 
$\Gamma_{\bet}$ is infinite. Then $X_{\bet}$ contains the Zariski closure of 
the corresponding coset, which is a translate of a positive dimension algebraic 
subgroup of $\mathbb{G}_a^g$. This would contradict the property inherited by 
$X_{\bet}$ from $X$. Thus $X_{\bet}(K_{\bet})\cap\Gamma_{\bet}^g$ is finite. 
Using \eqref{E:1900}, we conclude that $X(K)\cap\Gamma^g$ is finite.
\end{proof}

\begin{remark}
\label{R:Remark 1}  
Theorem \ref{T:T3} is an instance of Statement \ref{C:Con} 
because if we assume \eqref{C:Con} and we work with the hypothesis on $X$ 
from Theorem ~\ref{T:T3}, 
then, with the notations from \eqref{C:Con}, the intersection of $X$ 
with any translate of $B_i$ is finite. Otherwise, the Zariski 
closure of $X\cap (\gamma_i+B_i)$ would be a translate of a positive 
dimension algebraic subgroup of $\mathbb{G}_a^g$, and it would be contained in 
$X$.
\end{remark}

\end{document}